\documentclass[11pt]{amsart}
\usepackage[singlespacing]{setspace} % doublespacing

\usepackage{amsfonts}
\usepackage{mathrsfs} %\mathscr
\usepackage{dsfont}  %\matds

\usepackage[all]{xy}
\usepackage{stmaryrd}

\usepackage{bm,amsmath,amssymb,amsthm}
\usepackage{hyperref}
\theoremstyle{plain}
\newtheorem{theorem}{Theorem}

\swapnumbers
\newtheorem{proposition}[subsection]{Proposition}
\newtheorem*{proposition*}{Proposition}
\newtheorem{lemma}[subsection]{Lemma}
\newtheorem*{lemma*}{Lemma}
\newtheorem*{corollary*}{Corollary}
\swapnumbers

\theoremstyle{definition}

\newtheorem*{remark}{Remark}

\theoremstyle{remark}

\numberwithin{equation}{section}

\newcommand{\BmQ}{\mathbb{Q}}
\DeclareMathOperator{\MCl}{Cl}
\newcommand{\Mdemi}{\frac{1}{2}}
\newcommand{\Fmn}{\mathfrak{n}}
\newcommand{\SmA}{\mathscr{A}}
\newcommand{\CmL}{\mathcal{L}}
\newcommand{\BmZ}{\mathbb{Z}}
\newcommand{\SmC}{\mathscr{C}}
\newcommand{\BmR}{\mathbb{R}}
\newcommand{\SB}{\backslash}
\newcommand{\FmH}{\mathfrak{H}}
\newcommand{\CmO}{\mathcal{O}}
\newcommand{\Fma}{\mathfrak{a}}
\newcommand{\TmN}{\text{\bf N}}
\newcommand{\TmG}{\text{\bf G}}
\newcommand{\TmH}{\text{\bf H}}
\newcommand{\TmR}{\text{\bf R}}
\newcommand{\CmC}{\mathcal{C}}
\DeclareMathOperator{\Cl}{\TmC\Tml}
\newcommand{\TmC}{\text{\bf C}}
\newcommand{\Tml}{\text{\bf l}}
\newcommand{\Mdede}[4]{%
\begin{pmatrix}
#1&#2  \\
#3&#4 \\
\end{pmatrix}
}
\newcommand{\Mmod}[1]{\pmod{#1}}
\newcommand{\TmT}{\text{\bf T}}
\newcommand{\Fmb}{\mathfrak{b}}

\newcommand{\BmC}{\mathbb{C}}

\newcommand{\CmF}{\mathcal{F}}

\newcommand{\abs}[1]{\ensuremath{\left|#1\right|}}
\DeclareMathOperator{\MRe}{\Re e}
\DeclareMathOperator{\Mim}{\Im m}

\newcommand{\suma}{\sideset{}{_\SmA}\sum}
\providecommand{\remark}{\noindent {\bf Remark. }}

\begin{document}

\title{Heegner points and Eisenstein series} %le [] est pour ne pas avoir de running head

%\markleft{RIEN}
%------Information for first author-------%
\author[]{Nicolas Templier}   %le [] est pour ne pas avoir de running head
\address{D\'epartement de Mathematiques, Universit\'e Montpellier II, Place
  Eug\`ene Bataillon, 34095 Montpellier, FRANCE.}
\email{nicolas.templier@math.univ-montp2.fr}
%\curraddr{10, all\'ee des sciences appliqu\'ees, Res. R2 Apt. 1570,
%  31400 Toulouse, FRANCE.
%E-mail: nicolas.templier@math.univ-montp2.fr}
%\thanks will become a 1st page footnote.

%--------Information for second author------%
%\author[...]{...}
%\address{...}
%\curraddr{...}
%\email{...}
%\urladdr{...}
%-------- General info -------%
%\dedicatory{...}
\date{\today}
%\translator{...}
%-----------------------------%

\begin{abstract}
We give an alternative computation of the twisted second moment of critival values of class group $L$-functions attached to an imaginary quadratic field. The method avoids long calculations and yields the expected polynomial growth in the $s$-parameter for the remaining term.
\end{abstract}

\maketitle

\tableofcontents
%\newpage
%%%%%%%%%%%%%%%%%%%5\include{0Plan}
\section{Introduction.}

In this article, we are mainly concerned with the family of $L$-functions of characters on the class group of an imaginary quadratic field. We denote by the letter $D$ a fundamental negative discriminant in the sequel. We have the Kronecker symbol $\chi_D$, the quadratic field $K=\BmQ(\sqrt{D})$ and its associated ideal class group $\MCl_K=\MCl_D$ which is of size $h(D)$, the class number. We denote by $\chi$ a unitary character on this group (``class group character''), and $L(s,\chi)$ its associated Hecke $L$-series. 

The conductor of $L(s,\chi)$ is $|D|$ and the size of the family is $h(D)$ which is roughly $\abs{D}^{1/2}$. This makes the analysis of the moments a difficult task, so that only the first and second moments are understood for now. This knowledge is the key to deal with the nonvanishing problem, as shown by V.~Blomer \cite{Blom04}.

This family is interesting because of the simplicity of the geometric objects it is related with, namely Eisenstein series and Heegner points on modular surfaces, so that the computations are pleasant-looking. We know since Gauss that the class group is intimately related with definite binary quadratic forms. The roots of these forms (in the upper-half plane) are called CM or Heegner points. Moments of special values of $L$-functions associated to class group characters can often be interpreted as periods of particular automorphic form against Heegner points. In the present paper we focus on the special values $|L(s,\chi)|^2$ for $\MRe s=\Mdemi$, so that the period formula is well-known and due to Hecke. Our first result is the following:

%Our main observation \ref{maintool} is that when the growth
%in the cusp is very slow, the previous estimate still holds. The
%reason is that the highest in the cusp Heegner point is
%$\tau=\frac{1+\sqrt{D}}{2}.$ This enables us to deal with functions of
% logarithmic growth.

%investigate some averages of the values of
%modular functions on Heegner points. We consider here functions that
%are in $L^{1+\epsilon}$ for every $\epsilon>0$, but fails to be
%integrable because of the growth in the cusp for which precise
%estimate are known. We know at which rate the Heegner points climb in
%the cusp (see proposition 1), and comparing with this estimate, we
%find the (see proposition 2).

%Let's consider the function $\log |j(z)|$.

% we are able to apply spectral methods and
%Duke's equidistribution theorem, to obtain a sharp
%asymptotic with a power saving.

\begin{theorem}\label{th:dfiA}
There exist absolute constants $A,\delta>0$ such that we have uniformly in $D$ and $s$ with $\Re s=\Mdemi$:
\begin{multline*}
\frac{1}{h(D)}\sum\limits_{\chi\in\widehat{\MCl}_K}{\left|L(s,\chi)\right|^2}=\zeta(2)^{-1} L(1,\chi_D)\bigl[\CmL_D+\gamma-\log 2-2\frac{\zeta'}{\zeta}(2)+2\Re \frac{\xi'}{\xi}(2s)\bigr]|\zeta(2s)|^2 \\
    + \Re \frac{\Gamma(s)}{\Gamma(1-s)}(\frac{\sqrt{|D|}}{2\pi})^{2s-1}\zeta(2s)^3\zeta(4s)^{-1}L(2s,\chi_D)
+O(|s|^A|D|^{-\delta}),
\end{multline*}
where $\CmL_D:=\Mdemi \log |D| + \frac{L'}{L}(1,\chi_D)$, $\gamma$ is Euler constant and $\xi(s):=\pi^{-s/2}\Gamma(\frac{s}{2})\zeta(s)$ is the completed Riemann zeta function.
\end{theorem}

\remark It is not obvious that $\CmL_D$ is the leading term of the asymptotic. This is the case because, as we shall see in section~\ref{sec:LD}, we have $\CmL_D \gg \log |D|$. This lower bound is consequence of Burgess estimate.

\medskip 

Theorem~\ref{th:dfiA} is essentially due to W.~Duke, J.~Friedlander and H.~Iwaniec (Theorem $3$ from \cite{DFI4}), except that we have improved on two aspects. First, the remaining term is more precise here: {\em polynomial growth in the $s$-parameter}, whereas the growth was exponential in the initial article \cite{DFI4}. This improvement was expected, see the discussion at page 8 of \cite{DFI4}. In section \ref{discussion} we briefly discuss the importance of polynomial growth and state without proof some corollaries.

Second, our right-hand side is completely explicit. This yields closed formulas (a combination of special values and derivatives of $L$-functions) for the main term, whereas in \cite{DFI4} it is expressed in terms of intricate integrals. This main term becomes particularly complicate when specialized to $s=\Mdemi$. Namely, in the notations of~\cite{DFI4}:
\begin{equation}
\frac{1}{h(D)}
\sum\limits_{\chi\in\widehat{\MCl}_K}
\Bigl|L(\Mdemi,\chi)\Bigr|^2
=
 \sum_{0\le j+k\le 3} 
 c_{jk} L^{(j)}(1,\chi_D)(\log |D|)^k
 +O(|D|^{-\delta})
\end{equation}
From Theorem~\ref{th:dfiA}, we have explicit expression for the absolute constants $c_{jk}$ (directly in terms of derivatives of the $\Gamma$-function) if we observe that the right hand-side is real-analytic in the $s$ variable (on the line $\Re s=\Mdemi$). This must be so because both the remaining term and the left-hand side are regular on $\MRe s=\Mdemi$. One may also check by hands that the triple poles at $s=\Mdemi$ do cancel. It is interesting to observe that when $s=\Mdemi$ the main term is $L(1,\chi_D)(\log |D|)^3$, instead of $L(1,\chi_D)\log |D|$ when $s\not=\Mdemi$.

Previous aspect is of practical interest: when mollifying a family, an explicit expression for the main term of the asymptotic is required (as clearly explained in the discussion preceding Theorem 1.2 of \cite{KMV}). In the present case, Blomer \cite[Lemma 3.1]{Blom04} observed that a cancellation occurred in the explicit expression of $c_{00}$ given in \cite{DFI4}. At page \pageref{check} we check that the resulting expression is consistent with Theorem~\ref{th:dfiA}. Blomer further wrote that ``this lucky fact might indicate that there is a more elementary way of computing'' the main term in the asymptotic. Indeed the proof given in the present paper is short and avoids most unpleasant computations.

Our method is also convenient to deal with the twisted average. The following result with exponential instead of polynomial growth in the $s$-parameter, is due to Duke, Friedlander and Iwaniec (Theorem 4 from \cite{DFI4}).
\begin{theorem}\label{th:dfiB}
When $N$ a prime number and $D$ a fundamental discriminant with $\chi_D(N)=$~$1$ let us denote by $\Fmn$ an integral ideal such that $\Fmn\overline{\Fmn}=(N)$. There exist absolute constants $A,B,\delta>0$ such that the bound
\begin{equation*}
\frac{1}{h(D)}\sum\limits_{\chi\in\widehat{\MCl}_K}{\chi([\Fmn])}{\left|L(s,\chi)\right|^2} \ll_\epsilon
(N^{-\Mdemi}(\log N)^3 +N^A|D|^{-\delta}) |D|^\epsilon |s|^B
\end{equation*}
holds uniformly with respect to all the parameters ($N,D$ and $s$ with $\Re s=\Mdemi$).
\end{theorem}

\remark In some sense we are in presence of a ``bad family'' because the knowledge of Theorem~\ref{th:dfiB} yields subconvex bounds~\cite{DFI4} under conditional assumptions concerning the Landau-Siegel zero. A larger family is more suited to solve the subconvexity problem~\cite{DFI8}. This fact explains why researchers did not focus quite a lot on this family since the appearance of \cite{DFI4}.

\subsection{Discussion of the proofs} Our starting point is Hecke formula~\eqref{Hecke} which yields explicit formula~\eqref{hecke1} for the left-hand side in Theorem~\ref{th:dfiA}. It is then tempting to apply Duke's Theorem \cite{Duke88} which states that the Heegner points $(\tau^\SmA)_{\SmA\in\MCl_D}$ equidistribute on the modular surface when $D\rightarrow -\infty$. This idea, which we shall follow in the proofs, has been considered by Duke, Friedlander and Iwaniec themselves (see page 13 of \cite{DFI4}). The authors point out two problems:
\begin{quote} \it
``In the first place, the fact that the Eisenstein
series is not square-integrable causes technical difficulties. In the second place,
this method does not seem amenable to the twisted sums occurring in Theorem 4
and needed for the main applications. Therefore we shall use an alternative approach.''
\end{quote}

Indeed taking care of the singularities of the Eisenstein series is a rather tedious task. The function $z\mapsto |E(s,z)|^2$ is of moderate growth at infinity (in fact it is not $L^1$). Regularization process have been considered in the past and we briefly recall some of them in subsection~\ref{nearby}. Here we shall proceed by comparing the growth of a product of two unitary Eisenstein series and a linear combination of non-unitary Eisenstein series, see \eqref{diff1} and \eqref{diff2}. In subsection~\ref{subsec:Zagier}, we discuss in detail the difference with Zagier's regularization which is very close to the regularization carried out in the present paper.

After this regularization process is made, some singularities still remain. These are of logarithmic growth instead of polynomial growth. In section~\ref{sec:maintool} we state a simple proposition which is convenient to deal with these logarithmic singularities. The proof of Theorem~\ref{th:dfiA} will occupy section~\ref{sec:dfiA}.

The second difficulty alluded to by Duke-Friedlander-Iwaniec appears in the proof of Theorem~\ref{th:dfiB}. Both Heegner points $\tau^{\SmA}$ and $\tau^{\SmA[\Fmn]}$ appears in the right-hand side of formula \eqref{startn}: it seems difficult at first sight to apply equidistribution. We overcome this difficulty in section~\ref{sec:dfiB}. The right-hand side is very reminiscent of a recent work by V.~Vatsal, see e.g. \cite[Lemma 2.9]{Vatsal} We were inspired by the analytic part of this work, see \cite{Temp07} for an introduction: we shall study in detail the interplay between Heegner points of level $1$ and $N$ (in \cite{Vatsal} the interplay is between Heegner points of level $p^n$ and $p^{n+1}$). We show in section~\ref{subsec:sketch} that the addition of the extra twist $\chi([\Fmn])$ is equivalent to the action of the Hecke operator $\TmT_{\TmN}$.

\subsection{Nearby results in the literature}\label{nearby}
In \cite{MV}, Ph.~Michel and A.~Venkatesh observed that if $|L(s,\chi)|^2$ gets replaced by $L(\Mdemi,f\times \chi)$ where $f$ is a Maass form, equation \eqref{hecke1} still holds with $E(s,z)$ replaced by $f(z)$, a period formula recently established by S.-W.~Zhang \cite{Zhan01}. In that case the difficulty with the cusps is not present. It is thus possible to apply Duke's equidistribution Theorem directly. This applies in the same way for the family of ``canonical Hecke characters`` as recently shown by R.~Masri, see \cite{Masr07:quantitative}, \cite{Masr07:asymptotics} and the reference herein. In that case $E(s,z)$ gets replaced by $\theta_{d,k}(z)$, a (in general non-holomorphic) theta series.

In \cite{MY07}, Masri and T.~Yang deal with a function $\abs{\theta_{d,k}}^2$ which is $L^1$ but not of rapid decay (when $k$ is odd). Their proposition 6.1 resembles our proposition~\ref{maintool}. We explain briefly the differences in the remarks following proposition~\ref{maintool}.

In the theory of singular moduli initiated by R.~Borcherds and D.~Zagier, (exact!) evaluations of Heegner periods against meromorphic modular forms (or against weak Maass forms) have been discovered. A typical example is the trace of singular moduli:
\begin{equation}\label{tracemoduli}
 \sum_{\SmA\in \MCl_K} j(\tau^\SmA)
\end{equation}
where $j(z)=\dfrac{1}{q}+744+\cdots$ is the $j$-function. In \cite{Duke06}, Duke established a kind of asymptotic formula (\cite[(3)]{Duke06}) for this sum. He subtracted from $j$ a Poincar\'e series which is $1/q=e^{-2i\pi z}$ near $z=i\infty$ and applied afterwards the equidistribution theorem. Inspired by his result he also recovered Zagier's exact evaluation of \eqref{tracemoduli} via Fourier expansions of Poincar\'e series and Kloosterman/Sali\'e sums. This example of Duke is impressive because of the exponential growth of the $j$-function (to be compared to the moderate growth of the Eisenstein series $z\mapsto |E(s,z)|^2$). 

But usually how to solve a singularity problem heavily depends on the situation. And the regularization process from \cite{Duke06} is not well-suited to the present situation. Actually it is easy to obtain via a Poincar\'e series argument the asymptotic $\zeta(2)^{-1}L(1,\chi_D)\CmL_D+O(L(1,\chi_D))$ in Theorem~\ref{th:dfiA}, but not much better. Theorem 1.1 is stronger: it precises that the $O(1)$ is ''something`` $+O(|D|^{-\delta})$. This is the same discussion as in the introduction of \cite{DFI4}, where the authors points out that their Theorem 3 is much deeper than their Theorem 2. 

For these reasons we have developed another approach. Our approach is not very far away from section 3 of \cite{Zagi81}. Actually the closest paper related to the present work is \cite{TV84} (russian) by A.~Takhtajan and L.~Vinogradov.

\subsection{Comparison with Zagier's regularization}\label{subsec:Zagier}\footnote{We are grateful to a referee for pointing out to us the existence of that article.}
It is interesting to compare the content of our work with \cite{Zagi81}. We carry out the discussion with a lot of details because the comparison applies to other contexts. 
In \cite{Zagi81}, Zagier studies (and defines) the Rankin-Selberg transform
\begin{equation}\label{zagier}
 \int_{\Gamma\SB \FmH} F(z) E(s,z)\frac{dxdy}{y^2}
\end{equation}
for functions $F$ ``which are not of rapid decay``. For many choices of $F$ (confer the seven corollaries of the main Theorem in \cite{Zagi81}), this transform is equal (or closely related) to a Dirichlet $L$-series, as function in the $s$-variable. The main point is the analytic continuation, functional equation and poles of this Dirichlet $L$-series. These properties follow from the corresponding properties of the Eisenstein series.

One might understand the content of the present article in the following comparative way. We have ''replaced`` the Rankin-Selberg transform by periods over Heegner points:
\begin{equation}\label{intro:period}
\frac{1}{h(D)}\sum_{\SmA\in\MCl_K} F(\tau^\SmA).
\end{equation}
In our case, $F(z)=E(s_1,z)E(s_2,z)$ is the product of two unitary Eisenstein series (this is example 3 from \cite{Zagi81}). The period~\eqref{intro:period} is equal to the critical values studied in Theorems~\ref{th:dfiA} and \ref{th:dfiB}, see equations \eqref{hecke1} and \eqref{startn}. The main point here is the asymptotic behavior as $D\rightarrow -\infty$. The determination of this asymptotic behavior follows from Duke's equidistribution Theorem. 
 
We can pursue the analysis further. When $D$ gets large, \eqref{intro:period} approaches the ''mean`` of $F$ (this function $F$ is not $L^1$). We shall see (statement~\eqref{prop3} of proposition~\ref{prop:B} and statement \eqref{integrale}) that, after regularization process, this mean is zero. This is closely related to Example 2 of \cite{Zagi81} which reads:
\begin{equation}\label{meanzero}
 R.N.(\int_{\Gamma\SB \FmH} E(s_1,z)E(s_2,z) \frac{dxdy}{y^2})=0.
\end{equation}
The notation $R.N.$ is taken from \cite{Zagi81}. We shall not use this notation in the sequel because we define the regularized integral by a different approach. Equality~\eqref{meanzero} is a corollary of the main Theorem of \cite{Zagi81}. A cancellation occurs in the computation of this mean, and the fact that this mean is zero is important in the derivation of exact main terms in Theorems~\ref{th:dfiA} and \ref{th:dfiB}. As pointed out in section 2 of \cite{Zagi81}, equality \eqref{meanzero} may be derived from the Maass-Selberg relations. In fact we shall use the Maass-Selberg relations to establish statement~\eqref{prop3}.

However the bulk of the difficulty lies in a careful study of the error terms. We were inspired by a forthcoming article of Michel and Venkatesh~\cite{MV:GL2} where the authors (among many other things) give a clear framework for a {\em regularized Plancherel formula}. Here we would rather need a ''regularized spectral expansion'' of the function $F$. If this were available we could transform \eqref{intro:period} into (we display only the continuous part of the spectrum): 
\begin{equation}
 =\int_{\MRe s=\Mdemi} 
 \mu_P(ds)
 \cdot
 \left(\int_{\Gamma\SB \FmH} F(z) E(s,z)\frac{dxdy}{y^2}\right)
 \left(\frac{1}{h(D)}\sum_{\SmA\in\MCl_K} E(s,\tau^\SmA) \right)
 +
 \cdots
\end{equation}
where $\mu_P$ is the Plancherel measure. The first integral which is equal to \eqref{zagier} could be handled thanks to Zagier's method (see also section 5 of \cite{Zagi81} which builds a {\em regularized Petersson scalar product}). The period of the Eisenstein series (a Weyl's sum) is given by Hecke formula \eqref{hecke1}. In principle this would yield Theorem~\ref{th:dfiA}.

Unfortunately this ''regularized spectral expansion`` does not seem to exist in general. And it is clear that giving a rigorous framework to the process sketch above would be really lengthy task. For these reasons we cannot simply adapt the regularization of \cite{Zagi81}, \cite{Duke06} nor \cite{MV:GL2}. Instead we modify arguments in depth and produce independent proofs.

A last point of comparison is the following. We subtract from $F$ a main term, see equations \eqref{myreg1} and \eqref{regularized}. The point is that we know the exact evaluation of the period over Heegner points of the quantities we have subtracted. This turns out to be precisely the main term of the asymptotics. In equations (15) and (29) of \cite{Zagi81} a similar main term is subtracted. The main point is that the exact evaluation of the Rankin-Selberg transform of this quantity is known: it is zero, see Remark.~1 p.~420 in \cite{Zagi81}. 

%Of course this comparison is not surprising since the problems of analytic continuation of Dirichlet series %and of determining asymptotic behaviors are intimately related in analytic number theory.
% En fait la methode du present papier est tres proche de la formule de Plancherel regularisee. Plus precisement, on fait une decomposition spectrale regularisee.

\remark In Chapter 10 of \cite{thes:Temp} we studied the case where, in \eqref{intro:period}, the function $F$ is given by a Selberg kernel (this is related to section 7 of \cite{Zagi81}). This case is an application of Theorems~\ref{th:dfiA} and \ref{th:dfiB} as we shall explain briefly in section~\ref{discussion}; having polynomial versus exponential bounds is crucial there.

\subsection{Structure of the paper} In section~\ref{sec:notation}, we set up notations and recall classical formula of Hecke and Kronecker. In section~\ref{sec:LD} we recall fundamental subconvex bounds of Burgess and Duke-Iwaniec and establish $\CmL_D \gg \log |D|$. In section~\ref{sec:maintool} we are concerned with periods against Heegner points of functions with logarithmic singularities in the cusp. In sections~\ref{sec:dfiA} and \ref{sec:dfiB} we prove Theorems~\ref{th:dfiA} and \ref{th:dfiB} respectively.

\subsection{Acknowledgments} This article is based on Chapter 9 of the author's PhD thesis \cite{thes:Temp}. I thank Philippe Michel and Akshay Venkatesh for sharing with me some determinant ideas from \cite{MV:GL2} and for helpful discussions.
%%%%%%%%%%%%%%%%%%%%%%%%%%%%%%%%%%%%%%%%%%%%%%%%%%
\section{Hecke and Kronecker Formulas.}\label{sec:notation}
In this section we recall some classical facts and formulas about
Heegner points on $PSL_2(\BmZ)\SB \FmH$.
\subsection{Notations.} Set $\FmH$ for the upper half plane $\{x+iy,\ y>0\}$. The hyperbolic measure
$\dfrac{dxdy}{y^2}$ on $\FmH$ projects to $Y(1)$ or $Y_0(N)$ which are surfaces of finite volume. The normalized measure will be denoted $\mu(dz)$, so that $\mu(dz)=\dfrac{3}{\pi}\dfrac{dxdy}{y^2}$ (resp. $\dfrac{3}{\pi(N+1)}\dfrac{dxdy}{y^2}$) on $Y(1)$ (resp. $Y_0(N)$).

A clear reference for Heegner points is \cite[chap.~1]{GZ}. When $D\equiv 1 \pmod{4}$, the Heegner point which is the highest in the cusp is
\begin{equation*}
\tau:=\frac{-1+\sqrt{D}}{2}~;~\tau^2+\tau+\frac{D-1}{4}=0;
\end{equation*}
it is associated to the principal ideal class $[\CmO]$. The other Heegner points of discriminant $D$ are denoted $(\tau^\SmA)_{\SmA\in\MCl_K}$ to keep in mind the Galois action coming from the theory of Complex Multiplication. This is only a notation and we emphasis that we make nowhere use of results from that theory. We denote a sum over the class group by $\suma$.
%...........................%
\subsection{Hecke formula.} The value of an Eisenstein series $E(s,z)$ at a Heegner point is
\begin{equation}
\label{Hecke}
E(s,\tau^\SmA)=\frac{w}{2}(\frac{\sqrt{|D|}}{2})^{s}\cdot\zeta(2s)^{-1}\zeta(s,\SmA),
\end{equation}
where $w$ is half the number of units and which is $1$ when $|D|>4$, $\zeta$ is Riemann zeta function and $\zeta(s,\SmA)$ denote the partial zeta function associated to the ideal class $\SmA$
\begin{equation*}
 \zeta(s,\SmA):=\sum_{\Fma\in \SmA} (\TmN \Fma)^{-s}~;~
\suma \zeta(s,\SmA)=\zeta_K(s).
\end{equation*}
Here $\zeta_K(s)=\zeta(s)L(s,\chi_D)$ is the zeta function of $K=\BmQ(\sqrt{D})$.
%..........................%
\subsection{Kronecker limit formula.} Dedekind $\eta$-function and Ramanujan $\Delta$-function are related by the relation $\eta^{24}=\Delta$. The constant term of the Laurent expansion of the Eisenstein series at $s=1$ is known:
\begin{equation}\label{Kronecker}
\frac{\pi}{3}E(s,z)=\frac{1}{s-1}-\log|\Mim z\cdot\eta(z)^4|+2c+O(s-1)~;~
c:=\gamma-\log 2 -\frac{\zeta'}{\zeta}(2).
\end{equation}
\subsection{An exact average.} As consequence of the two previous formulas we infer that
\begin{equation}\label{Main1}
\frac{-1}{h(D)}\suma {\log |\Mim \tau^\SmA
 \cdot \eta(\tau^\SmA)^4|}=\CmL_D+\log 2-\gamma.
\end{equation}
We think of this formula as a reference that quantifies the climb in the cusp of the Heegner points: the moderate growth of $z\mapsto \log |\Mim z\cdot \eta(z)^4|$ in the cusp makes it not $L^1$ but barely ; its average over Heegner points is a quantity which grows with $D$ at logarithmic rate.
%%%%%%%%%%%%%%%%%%%%%%%%%5\include{7burgess}
\section{Subconvex Bounds.}\label{sec:LD}
In this section we recall the fundamental subconvex estimates of Burgess and Duke-Iwaniec. We also study in detail the quantity $\CmL_D$.
\subsection{Burgess bound.}
\begin{theorem}[Burgess]\label{Burgess}
\begin{align}\label{Ba}\tag{B1}
L(s,\chi_D) \ll_\epsilon |D|^{\tfrac{3}{16}+\epsilon}\cdot|s|,&\quad\text{for }\Re s=\Mdemi.\\
\label{Bb}\tag{B2}
\sum\limits_{n=M}^{N+M}{\chi_D(n)} \ll_{\epsilon,r} N^{1-\frac{1}{r}}|D|^{\tfrac{r+1}{4r^2}+\epsilon},&\quad\text{for any integer $r\ge 1$.}
\end{align}
\end{theorem}

The quantity $\CmL_D$ has appeared before in the literature in various contexts. We exhibit an unconditional lower bound \eqref{log} which has the right order of magnitude. This is equivalent to have the $\liminf$ in \eqref{eq:prop:liminf} greater than $-\Mdemi$. That the $\liminf$ is greater or equal to $-\Mdemi$ is a consequence of the Polya-Vinogradov inequality (see \cite[p.~326]{Iwan04}). To go beyond, it is necessary to use Burgess subconvex estimate, and in  \cite{Mich01}, $-\frac{3}{8}$ is obtained thanks to~\eqref{Ba}. The best result we display below (the bound $-\frac{1}{4}$) is taken from a remark in unpublished notes of Iwaniec and uses the full Burgess estimate~\eqref{Bb}.
\begin{proposition} We have the following
\begin{equation}\label{eq:prop:liminf}
\liminf_{|D|\rightarrow \infty} \frac{L'(1,\chi_D)}{\log |D|\cdot
  L(1,\chi_D)} \ge -\frac{1}{4}.
\end{equation}
\end{proposition}
A corollary is that
\begin{equation}\label{log}
\CmL_D \ge (\frac{1}{4}-\epsilon)\log |D|
\end{equation}
for any $\epsilon>0$ and $|D|$ large enough (depending on $\epsilon$).
Hence $\CmL_D$ is indeed the leading term in Theorem~\ref{th:dfiA} and in formula~\eqref{Main1}.

\remark Under $\TmG\TmR\TmH$:
\begin{gather*}
\CmL_D=\Mdemi \log |D| + O(\log\log |D|)~;~
\lim_{|D|\rightarrow \infty} \frac{L'(1,\chi_D)}{\log |D|\cdot
  L(1,\chi_D)}=0.
\end{gather*}
%\begin{proof} Adapted from \cite{Mich01}.
%\begin{equation}
%\sum\limits^\infty_{n=1}{\frac{r(n)}{n}V(\frac{n}{X})}=\int_{(\sigma)}{\zeta(s+1)L(s+1,\chi_D)\cdot\widehat{V}(s)X^s\frac{ds}{2i\pi}},
%\end{equation}
%where $\sigma>1$ ; $r(n)=\sum\limits_{ab=n}{\chi_D(a)}$ is the number of ways of writing the integer %$n$ as the norm of an element in $K$ ; $V\in\CmC^\infty_c(\BmR_+)$ is a fixed test function, non negative and which takes the value $1$ near $0$ ; and $\widehat{V}$ is its Mellin transform,
%\begin{equation*}
%\widehat{V}(s)=\int_0^\infty{y^{s-1}V(y)dy}~;~
%V(y)=\int_{(3)}{\widehat{V}(s)y^{-s}\frac{ds}{2i\pi}}.
%\end{equation*}
%Moving the contour to the left we catch a pole at $s=1$, whose residue is
%\begin{equation}
%L(1,\chi_D)\log(\eta X) + L'(1,\chi_D)
%\end{equation}
%and the remainding term is
%\begin{equation}\label{intdemi}
%\int_{(\Mdemi)}{\zeta(s)L(s,\chi_D)\cdot\widehat{V}(s-1)\cdot X^{s-1}\cdot\frac{ds}{2i\pi}} %\ll_\epsilon |D|^{\frac{3}{16}+\epsilon} X^{-\Mdemi},
%\end{equation}
%applying \eqref{Ba}. We take $X=|D|^{\frac{3}{8}-\epsilon}$ and notice that the initial sum is %nonnegative - \cite{Mich01,Iwan04} prove that it is bounded below by a constant, but this is useless.

\begin{proof}
Our source is a remark in unpublished notes of Iwaniec that did not appear in \cite{Iwan04}.

\noindent We start with the sum:
\begin{equation}
\sum\limits_{n\le X}{\frac{r(n)}{n}}=\sum\limits_{a\le Y}{\frac{\chi_D(a)}{a}}\sum\limits_{b\le \frac{X}{a}}\frac{1}{b} ~+~ \sum\limits_{Y<a,~ab\le X}{\frac{\chi_D(a)}{ab}}
\end{equation}
that we have cut into two parts. The first one is:
\begin{equation*}
L(1,\chi_D)\cdot (\log X + \gamma) + L'(1,\chi_D)+O(\frac{Y}{X}+\frac{|D|^{\Mdemi+\epsilon}}{Y})
\end{equation*}
and the second one is bounded by:
\begin{equation*}
\ll_{\epsilon,r} |D|^{\tfrac{r+1}{4r^2}+\epsilon}X^{1-\tfrac{1}{r}+\epsilon}
\end{equation*}
which follows by partial summation and from \eqref{Bb}.

We choose $Y=X^{1-\epsilon'}$ and $X=|D|^{\tfrac{r(r+1)}{4r^2}+\epsilon''}$ and we obtain that for $|D|$ large enough,
\begin{equation*}
(\frac{r(r+1)}{4r^2}+\epsilon''/2)\cdot L(1,\chi_D)\cdot\log(|D|) + L'(1,\chi_D) \ge 0
\end{equation*}
The case $r=1$ should be thought as Polya-Vinogradov inequality. The case $r=2$ should be thought as \eqref{Ba} and yields the lower bound $-\frac{3}{8}$ from \cite{Mich01}. The good choice which concludes the proof of the Lemma is to take $r$ arbitrary large.
\end{proof}
%\remark One should think the second proof as a better way for giving a bound for \eqref{intdemi}. The $s$-integral enables to use the best Burgess estimates \eqref{Bb} whereas the subconvex bound \eqref{Ba} only make use of $r=2$.
\subsection{Duke-Iwaniec bounds.} We recall the fundamental Theorem obtained in \cite{Duke88,Iwan87} (the polynomial growth in $N$ appears in \cite[sections 13-14]{DFI4}).
\begin{theorem}\label{Duke}
The Heegner points become equidistributed on $Y_0(N)$ as $D\rightarrow -\infty$. More precisely, the Weyl's sum are bounded by
\begin{align*}
\frac{1}{h}\suma{E_\Fma(s,\tau^\SmA)}&\ll |s|^A N^B |D|^{-\delta};\\
\frac{1}{h}\suma{\phi_j(\tau^\SmA)}  &\ll t_j^A N^B |D|^{-\delta}
\end{align*}
for some absolute constants $A,B, \delta>0$ ; we have $\Re s=\Mdemi$ and $E_\Fma(s,z)$ is the Eisenstein series associated to the cusp $\Fma$ ; $\{\phi_j\}_j$ is an orthonormal basis of Maass cusp forms with respective Laplace eigenvalue $\dfrac{1}{4}+t^2_j$.
\end{theorem}

\section{Weyl's Sums with Logarithmic Singularities.}\label{sec:maintool}
An essential ingredient we shall make use of is the following estimate for Weyl's sum
of test function that are not perfectly smooth. In fact we allow some logarithmic
singularities at points which are slowly approached by the Heegner
points, typically a cusp.

\begin{proposition}\label{maintool}
Let $A(z)$ be a function on $Y(1)$ that satisfies the following assumptions. For $z=x+iy$ in the standard fundamental domain $\CmF$ of $PSL_2(\BmZ)$,
\begin{gather}\label{A1}\tag{A1}
\text{$A$ is of class $\SmC^\infty$ and }
\frac{\partial^{i+j}}{\partial x^i\partial y^j}{A(z)}\ll_{ij} y^{A_{ij}}\quad \text{for some $A_{ij}>0$};\\
\label{A2}\tag{A2}
A(z)\ll \max(\log y,1);\\
\label{A3}\tag{A3}
\int_\CmF A(z) \mu(dz)=0.
\end{gather}
There exists an absolute constant $\delta>0$, such that
\begin{equation}\label{eq:prop:A}
\frac{1}{h(D)}\suma {A(\tau^\SmA)}\ll D^{-\delta}.
\end{equation}
Moreover the constant involved in the last bound depends linearly
on the implicit constant appearing in \eqref{A2} and a fixed and finite number of implicit constants appearing in \eqref{A1}.
\end{proposition}

\remark We can also add a finite number of logarithmic singularities at
some fixed CM points. Indeed they are approached at $\log |D|$ rate by
the Heegner points of discriminant $D$, like the cusp. The proof is the same but with incomplete Eisenstein
series now constructed from a function $\psi$ centered at the height of the singularity.

\remark We have stated the theorem for the full modular surface $Y(1)$ but the proof works also on $Y_0(N)$ with obvious modifications.

\remark Proposition~\ref{maintool} resembles proposition 6.1 of \cite{MY07}. There are some differences. In \cite{MY07}, the growth is in $y^{1/2}$ which is more general that~\eqref{A2}. On the other hand our assumption~\eqref{A1} is more robust than the assumption ``$\Delta^a(A(z) - y^{1/2})$ is of exponential decay for any $a$'' in \cite{MY07}. The point in proposition 6.1 of \cite{MY07} is that the quantity subtracted $y^{1/2}$ does not depend on the $x$ variable, so that it is annihilated by $y^2\frac{\partial^2}{\partial x^2}$, a part of the Laplacian that introduces in general annoying growth at infinity. Also in our proof the truncation at height $Y:=|D|^\eta$ with $\eta$ very small whereas the truncation is at height $|D|^{1/2}$ in \cite{MY07}.

\begin{proof}
 Fix once and for all a nonnegative function $\psi_0\in\CmC^\infty(\BmR_+)$ with support on  $[1,+\infty[$ that takes the value $1$ in a neighborhood of infinity ; for $Y>0$, let $\psi(y):=\psi_0(\frac{y}{Y})$. All the constructions below are linear functional in $A$ and thus the constants involved in the various bounds depend linearly on the constants appearing in \eqref{A1} and \eqref{A2}.

Let's form the incomplete Eisenstein series (see \cite[p.~62]{Iwaniec}):
\begin{equation}\label{Mellin}
E(z|\psi):=\sum\limits_{\gamma\in\Gamma_\infty\SB \Gamma}{\psi(\Mim
      \gamma z)}
=\int_{(3)}{E(v,z)\widehat{\psi}(v)\frac{dv}{2i\pi}}
=\widehat{\psi}(1)+\int_{(\sigma)}{E(v,z)\widehat{\psi}(v)\frac{dv}{2i\pi}},
\end{equation}
for any $\sigma$ in $(0,1)$ and where we have introduced the Mellin transform,
\begin{equation*}
\widehat{\psi}(v):=\int_{0}^{\infty}{\psi(y)y^{-v-1}dy}.
\end{equation*}
We use it to isolate the contribution from the cusp, writing
\begin{equation}
A(z)=\bigl[1- E(z|\psi)\bigr]A(z)+ E(z|\psi)A(z)
\end{equation}
and then cut the average \eqref{eq:prop:A} accordingly as $S_1+S_2$.
Integrating by parts, we infer that
\begin{equation*}
\widehat{\psi}(v)\ll_{\sigma,A} |v|^{-A}Y^{A}
\end{equation*}
for all $A\ge 0$ and $\Re v=\sigma.$
Plugging the Burgess bound \eqref{Ba} in \eqref{Mellin}, we obtain
\begin{equation*}
\frac{1}{h}\suma{E(\tau^\SmA|\psi)}\ll Y^{-1}+|D|^{-\delta}Y^B
\end{equation*}
for some $B$, because $\widehat{\psi}(1)\ll Y^{-1}$.
Observe now that $E(\cdot|\psi)$ takes nonnegative values. Hence \eqref{A2} implies
\begin{equation}
S_2\ll \log |D|\cdot Y^{-1}+|D|^{-\delta}Y^B.
\end{equation}

The first sum $S_1$ is addressed through equidistribution. Since we need a quantitative bound, we have to make our
Weyl sums explicit by expanding spectrally $A(\cdot)(1- E(\cdot|\psi))$. With \eqref{A1}, we can perform successive integrations by part to bound the $r$-spectral coefficient by
$\ll r^{-A}Y^{B}$. Then the Duke-Iwaniec estimates (Theorem~\ref{Duke}) show that the total contribution of the non trivial Weyl sums is bounded by $\ll |D|^{-\delta} Y^B$.

By \eqref{A2}, $\int A(z)E(z|\psi)\mu(dz) \ll Y^{-1}$.  And this together with \eqref{A3}
provides a bound for the principal Weyl sum and implies
\begin{equation}
S_1 \ll Y^{-1}+|D|^{-\delta}Y^B.
\end{equation}

The choice $Y=|D|^\eta$ for $\eta>0$ sufficiently small achieves the
proof of Proposition~\ref{maintool}.
\end{proof}

%%%%%%%%%%%%%%%%%%5\include{5theorem2}
\section{Proof of Theorem~\ref{th:dfiA}.}\label{sec:dfiA}
By \eqref{Hecke} and orthogonality of characters we infer that
\begin{equation}\label{hecke1}
\frac{1}{h}\sum\limits_{\SmA\in\MCl_K}{|\zeta(2s)E(s,\tau^\SmA)|^2}=w^{2}(\frac{\sqrt{|D|}}{2})\frac{1}{h^2}{\sum\limits_{\chi\in
    \widehat{\Cl}_K}}{|L(s,\chi)|^2}.
\end{equation}
We shall now estimate the left hand-side, building on standard properties of the Eisenstein series.
%\begin{align}\label{carres}
%\frac{1}{h}\suma{|\zeta(2s)E(s,\tau^\SmA)|^2}=&\frac{6}{\pi}(\CmL_D-\log 2 -
%2\frac{\zeta'}{\zeta}(2)+\gamma)\\
%&+\Re F_D(s)+\frac{12}{\pi}\Re
%\frac{\xi'}{\xi}(2s)+O(|s|^A|D|^{-\delta})\nonumber
%\end{align}
%where
%\begin{gather*}
%F_D(s):=\pi\frac{\Gamma(s)}{\Gamma(1-s)}\left(\frac{\sqrt{|D|}}{2\pi}\right)^{2s-1}\zeta(2s)^2\zeta(2-2s)^{-1}\zeta(4s)^{-1}L(1,\chi_D)^{-1}L(2s,\chi_D).
%\end{gather*}
%\remark The error term is poor when $s$ approaches $\Mdemi$. Actually a variation of the proof would yield the asymptotic of $\sum_\chi{|L(\Mdemi,\chi)|^2}$, and the main term contains a further $\log |D|$. We shall not do it here since this has already been obtain in \cite{DFI4}.

%\remark Observe that for $\Re s=\Mdemi$,
%\begin{gather*}
%\xi(2-2s)=\xi(2s-1)=\overline{\xi(2s)};\\
%\frac{\xi'}{\xi}(2-2s))=-\frac{\xi'}{\xi}(2s-1)=\overline{\frac{\xi'}{\xi}(2s)};\\
%F_D(s)=-F_D(\Mdemi-s)=\overline{F_D(1-s)}.
%\end{gather*}

Let's introduce for $\Re s=\Mdemi$, and $z\in Y(1)$ the function $B(s,z)$ defined by
\begin{equation}\begin{aligned}\label{myreg1}
B(s,z)&:=|\zeta(2s)E(s,z)|^2-2\Re \pi^{1-2s}\frac{\Gamma(s)}{\Gamma(1-s)}\zeta(2s)^2E(2s,z)-\\
&\qquad -\frac{6}{\pi}|\zeta(2s)|^2\bigl[-\log |\Mim z\cdot \eta(z)^4|+2c \bigr]-\frac{12}{\pi}|\zeta(2s)|^2\Re\frac{\xi'}{\xi}(2s).
  \end{aligned}\end{equation}
%\begin{equation}
%A(s,z):=|E(s,z)|^2 - \frac{6}{\pi}[-\log|\Mim z\cdot\eta(z)^4|+2c]
%-\Re \frac{\xi(2s)}{\xi(2-2s)}E(2s,z)-\frac{12}{\pi}\Re \frac{\xi'}{\xi}(2s)
%\end{equation}
This function is analytic in $s$ and it is possible to check that it has no singularity at $s=\Mdemi$ despite the presence of the poles of zeta (see below). Now Theorem~\ref{th:dfiA} follows from \eqref{Hecke}, \eqref{Main1}, \eqref{hecke1} and the following claim:
\begin{equation}\label{TermeErreur}
\frac{1}{h}\suma{B(s,\tau^\SmA)}=O(|s|^A|D|^{-\delta}).
\end{equation}
This in turn will be a consequence of Proposition~\ref{maintool} and so need to prove:
\begin{proposition}\label{prop:B} The function $B$ satisfies the following uniform estimates.
\begin{multline}\label{prop1}\tag{B1}
\qquad \quad
\text{$B(s,\cdot)$ is of class $\SmC^\infty$ for all $s$, and}\\
\frac{\partial^{i+j}}{\partial x^i\partial y^j}{B(s,z)}
\ll_{ij} 
|s|^{A_{ij}}\cdot y^{B_{ij}}\quad \text{for some $A_{ij},B_{ij} >0$};
\qquad
\end{multline}
\begin{equation}\label{prop2}\tag{B2}
B(s,z)\ll |s|^A\cdot\max(\log y,1)\quad \text{for some $A>0$};
\end{equation}
\begin{equation}\label{prop3}\tag{B3}
\int_\CmF B(s,z) \mu(dz)=0\quad \text{for all $s$}.
\end{equation}
\end{proposition}
\begin{proof}
 For the second assertion \eqref{prop2} we need to know the behavior of the Eisenstein series in the cusp which is recalled in the Lemma below. For the other claims, we introduce the function
 \begin{multline}\label{diff1}
  B(s_1,s_2,z):=\zeta(2s_1)\zeta(2s_2)E(s_1+s_2,z)+\pi^{2s_2-1}\frac{\Gamma(1-s_2)}{\Gamma(s_2)}\zeta(2s_1)\zeta(2-2s_2) E(1+s_1-s_2,z)\\
  +\pi^{2s_1-1}\frac{\Gamma(1-s_1)}{\Gamma(s_1)}\zeta(2-2s_1)\zeta(2s_2) E(1+s_2-s_1,z)\\
  +\pi^{2s_1+2s_2-2}\frac{\Gamma(1-s_1)\Gamma(1-s_2)}{\Gamma(s_1)\Gamma(s_2)}\zeta(2-2s_1)\zeta(2-2s_2) E(2-s_1-s_2,z)
\end{multline}
 and we write $B(s,z)$ as the 'value' at $s_1=s$, $s_2=1-s$ of
\begin{equation}\label{pf:regul}
   \zeta(2s_1)\zeta(2s_2)E(s_1,z)E(s_2,z)-B(s_1,s_2,z)
 \end{equation}
 We are allowed to do so because these functions are actually holomorphic in the strip $\frac{1}{4}< \Re s_1, \Re s_2 <\frac{3}{4}$ as a consequence of the Lemma below (it is also possible and probably better to make use of the symmetries $s_1\leftrightarrow 1-s_1$ and $s_2\leftrightarrow 1-s_2$ of $B(s_1,s_2,z)$ and the fact that the potential poles are at most simple). Let's go back to assertion \eqref{prop1}, put $s=\Mdemi+it$ and write the 'value' with the help of a Cauchy integral in $s_2$ on the rectangle with edges $\frac{5}{8}-(|t|+1)i$ ; $\frac{5}{8}+(|t|+1)i$ ; $\frac{3}{8}+(|t|+1)i$ ; $\frac{3}{8}-(|t|+1)i$. After this, we can bound each term of \eqref{diff1} individually, making use of convexity bounds for the entire function $s(1-s)\zeta(2s)E(s,z)$ (thanks to the Cauchy path, we are at bounded distance of any pole).

 For assertion \eqref{prop3} we take two generic $s_1$ and $s_2$ (that is with $s_1\not= s_2,1-s_2$ and $s_1,s_2\not=\Mdemi$) and we prove that the integral of \eqref{pf:regul} is zero ; then we make $s_1$ and $s_2$ tends to $s$, $1-s$ and use continuity (we could have made use of the same Cauchy integral instead). The function in \eqref{pf:regul} is $L^1$ so that we may introduce a cut-off at height $Y$ and let $Y$ tends to infinity ; the corresponding integral of each term of \eqref{diff1} and \eqref{pf:regul} can be computed explicitly thanks to the Maass-Selberg relations (see \cite[Proposition 6.8]{Iwaniec}) ; all the leading terms cancel out and the limit as $Y\rightarrow \infty$ is indeed zero \footnote{The paper \cite{MV:GL2} builds a very general framework which explains better why the integral \eqref{prop3} does vanish. We are grateful to the authors for bringing it to our attention.}.
\end{proof}
\begin{lemma}\label{lem:eis} The function in the $s$-variable
\begin{equation}\label{function}
  \zeta(2s)E(s,z)-\zeta(2s) y^s - \pi^{2s-1}\frac{\Gamma(1-s)}{\Gamma(s)}\zeta(2-2s)y^{1-s}
\end{equation}
has a holomorphic continuation to the whole complex plane, and satisfies for all positive integer $M$ the bound
\begin{equation}\label{lem:bound}
  \ll_{M,\epsilon} |s|^{M+1-\Re s}y^{-M+\epsilon}
\end{equation}
uniform on $-M\le \Re s \le M+1$ and $z=x+iy\in \CmF$ (fundamental domain).
\end{lemma}
We couldn't locate this Lemma in the literature where one usually writes down the exponential decay in the $z$-variable which is valid for any $s$, but with a constant growing exponentially in $s$. Here we have a weaker decay in the $z$-variable but a polynomial growth in $s$ which is precisely what is needed for Theorem~\ref{th:dfiA}. The result is certainly well-known, but we provide a quick proof in view of the importance of this lemma for Theorem~\ref{th:dfiA}.

\begin{proof}
  We know already that the function~\eqref{function} has a meromorphic continuation to the whole complex plane with potential poles of order at most $1$ at $s=0,\Mdemi,1$. One can check by various means that the residues actually cancel out. 
  
  For the second claim, we establish easily the bound~\eqref{lem:bound} for $\Re s=M+1$. On this line we appeal to the definition of the Eisenstein series \eqref{defeis} and we plug the standard estimates which establish its absolute convergence, see for instance \cite[Lemma 2.10]{Iwaniec} and \cite[(8.10-8.11)]{Ven05}.
  
  Then by the functional equation of the Eisenstein series, the bound~\eqref{lem:bound} is also valid for $\MRe s=-M$. 
  
  We conclude with the Phragm\`en-Lindel\"of principle which imply (by holomorphy) that the bound~\eqref{lem:bound} is valid on the whole strip $-M\le \Re s \le M+1$.
  \end{proof}

%-----------check--------------------------
\subsection{Comparison with \cite{DFI4}.} \label{check} 
We check that the asymptotic stated in Theorem 3 of \cite{DFI4} indeed agrees with our Theorem~\ref{th:dfiA}. We adopt the notations and labelling from \cite{DFI4} in this paragraph only. The main term is denoted there by $l_D(s)+l_D^+(s)$ ; the term $l_D(s)$ is given by (1.16) where the constants $c_{ij}(s)$ are extracted from (1.18), (6.26) and (6.24) ;  the off-diagonal term $l_D^+(s)$ is given in (11.17). Observe first that $c_{10}(s)=2c_{01}(s)$, so that the leading component is indeed governed by the quantity $\CmL_D=\frac{1}{2}\log |D|+\frac{L'}{L}(1,\chi_D)$. Next the coefficient $c(s)$ corresponds to our last term with product of zeta functions. The last coefficient $c_{00}(s)$ which is extracted from (6.24) is the longest to compute:
\begin{equation*}
c_{00}(s)=\frac{12}{\pi w}|\Gamma(s)|^{-2}\Bigl[2\Theta(s)\Theta(1-s)(-\log 2 -2 \frac{\zeta'}{\zeta}(2)+\gamma)+\Theta(s)\Theta'(1-s)+\Theta'(s)\Theta(1-s)+\widetilde{R}(1)\Bigr],
\end{equation*}
where $\Theta(s)=\xi(2s)$ in our notations and $\widetilde{R}(1)$ is an integral given in (6.20) and (6.21) which we don't recall here. The off-diagonal term $l_D^+(s)$ is given by (1.17) where another integral $I(s)$ shows up. It turns out that
\begin{equation*}
I(s)=-\frac{\pi}{12s^2(1-s)^2}=-\tilde{R}(1)
\end{equation*}
by moving the integral and computing a residue, so that these two last terms which did not appear in our Theorem~\ref{th:dfiA} indeed cancel out.

%%%%%%%%%%%%%%%%%%%%%%%%5\include{6theorem3}
\section{Proof of Theorem~\ref{th:dfiB}.}\label{sec:dfiB}

We keep the notations of the Theorem: $N$ is a prime number and $\Fmn$ is one of the two ideals such that $\Fmn\overline{\Fmn}=(N)$. Then we infer, by orthogonality of characters:
\begin{equation}\label{startn}
w^{2}(\frac{\sqrt{|D|}}{2})\frac{1}{h^2}|\zeta(2s)|^{-2}\sum\limits_{\chi\in\widehat{\MCl}_K}{\chi([\Fmn])\left|L(s,\chi)\right|^2}=\frac{1}{h}\sum\limits_{\SmA\in\MCl_K}{E(s,\tau^\SmA)\overline{E(s,\tau^{\SmA[\Fmn]})}}.
\end{equation}

\subsection{Modular curves - notations.}
The modular groups are $\Gamma:=\Gamma(1):=PSL_2(\BmZ)$ ;
$\Gamma_0(N):=\{\Mdede{a}{b}{c}{d},~c\equiv 0\Mmod{N}\}$ ; put $Y(1)=\Gamma\backslash\FmH$ ; the modular surface $Y_0(N)=\Gamma_0(N)\backslash \FmH$ has two cusps $0$ and $\infty$.
We have the important Atkin-Lehner-Fricke involution
$\sigma:=\Mdede{0}{-1}{N}{0}$. It is associated with the cusp $0$ in the sense that it conjugates the stabilizer of $0$, $\Gamma_0(N)_0$ with $\Gamma_\infty=\Mdede{1}{\BmZ}{0}{1}$. Observe also that
\begin{equation}\label{sigma}
\sigma^{-1}\Gamma\sigma\cap \Gamma=\Gamma_0(N).
\end{equation}
We introduce the two maps
\begin{equation*}
\xymatrix{
    & Y_0(N) \ar[ld] \ar[rd] &\\
    Y(1) & & Y(1)
  }
\end{equation*}
where the first one is the natural projection while the second one is
\begin{equation}
z\mapsto Nz\equiv \sigma z \Mmod{\Gamma}.
\end{equation}
We denote by $\iota:Y_0(N)\rightarrow Y(1)\times Y(1)$ the embedding induced by these two maps. This is the graph of the Hecke correspondence $\TmT_N$,
which is of degree $N+1=[\Gamma(1):\Gamma_0(N)]$.
\subsection{Eisenstein series.}
Recall that for $\Re s>1$
\begin{equation}\label{defeis}
E(s,z):=\sum_{\gamma\in \Gamma_\infty\backslash \Gamma} (\Mim
\gamma z)^s
\end{equation}
is the Eisenstein series on $Y(1)$ and that
\begin{equation}
E_\infty(s,z):=\sum_{\gamma\in \Gamma_\infty\backslash \Gamma_0(N)} (\Mim
\gamma z)^s~;~
E_0(s,z):=\sum_{\gamma\in \Gamma_0(N)_0\backslash \Gamma_0(N)} (\Mim
\sigma\gamma z)^s
\end{equation}
are the two Eisenstein series on $Y_0(N)$ associated to the cusp $0$ and $\infty$
respectively.
 They are intertwined by the following simple relations, all consequences of \eqref{sigma}.
\begin{equation}\begin{split}\label{eisensteinniveaux}
E(s,z)&=E_\infty(s,z)+N^{s}E_0(s,z);\\
E_\infty(s,z)&=(N^{s}-N^{-s})^{-1}\bigl[E(s,Nz)-N^{-s}E(s,z)\bigr];\\
E_0(s,z)&=E_\infty(s,\frac{-1}{Nz})=(N^{s}-N^{-s})^{-1}\bigl[E(s,z)-N^{-s}E(s,Nz)\bigr]\\
\end{split}\end{equation}

 It is possible to prove with Lemma~\ref{lem:eis} that the difference
\begin{equation}\begin{gathered}\label{growth}
\zeta^{(N)}(2s)E_{\Fma}(s,\sigma_\Fmb z)-\zeta^{(N)}(2s)\delta_{\Fma\Fmb}y^s-\zeta^{(N)}(2s)\phi_{\Fma\Fmb}(s)y^{1-s}
\end{gathered}\end{equation}
is entire in the $s$-variable and bounded by $\ll_M (N|s|)^{M+1} y^{-M}$. Here the scattering coefficients are given by
\begin{equation*}
\Mdede{\phi_{\infty\infty}(s)}{\phi_{\infty0}(s)}{\phi_{0\infty}(s)}{\phi_{00}(s)}:=\frac{\xi(2s-1)}{\xi(2s)}(N^{2s}-1)^{-1}\Mdede{N-1}{N^s-N^{1-s}}{N^s-N^{1-s}}{N-1}=M(s)^{-1}M(1-s)
\end{equation*}
where
\begin{equation*}
M(s):=\xi(2s)\Mdede{1}{N^s}{N^s}{1}.
\end{equation*}

Recall that $E_\Fma(s,z)$ has a simple pole at $s=1$ with residue $\dfrac{3}{(N+1)\pi}$.
%%%%%%%%%%%%%%%%%%%%%%%%%%%%%%%%%%
\subsection{Heegner points.}
We borrow facts from \cite[chap.~I]{GZ}.
There are $h(D)$ Heegner points,
$(\tau^{\SmA})_{\SmA\in\MCl_K}$ on the modular curve $Y(1)$.
There are $2h(D)$ Heegner points, $(\tau_\Fmn^{\SmA})$ on
the modular curve $Y_0(N)$, where ${\SmA\in\MCl_K}$ and $\Fmn\overline{\Fmn}=(N)$.
We have a modular interpretation as CM elliptic curves, respectively as cyclic isogenies between CM elliptic curves:
\begin{equation*}
\tau^{[\Fma]}=\BmC/\Fma~;~
\tau_\Fmn^{[\Fma]}=(\BmC/\Fma \rightarrow \BmC/\Fma\Fmn^{-1}).
\end{equation*}

We may also write down-to-earth formulas with coordinates:
to a point $[z]\in Y(1)$ corresponds the elliptic curve
$\BmC/<z,1>$ ; and to a point $[z]\in Y_0(N)$ corresponds the diagram of
elliptic curves $(\BmC/<z,1>\rightarrow \BmC/<z,\frac{1}{N}>)$.

On the modular curve $Y(1)$ we have:
\begin{equation*}
\tau^{[\Fma]}=\frac{B+\sqrt{D}}{2A} \Mmod{\Gamma}~;~
\Fma=A\BmZ+\frac{B+\sqrt{D}}{2}\BmZ~;~
\TmN \Fma=A.
\end{equation*}
On $Y_0(N)$, we have the following formulas:
\begin{gather*}\Fmn=N\BmZ+\frac{\beta+\sqrt{D}}{2}\BmZ~;~
\beta^2\equiv D \Mmod{4N}.\\
\tau_\Fmn^{[\Fma]}=\frac{B+\sqrt{D}}{2A} \Mmod{\Gamma_0(N)}~;~
\Fma=A\BmZ+\frac{B+\sqrt{D}}{2}\BmZ\\
\TmN \Fma=A\equiv 0\Mmod{N},~~B\equiv \beta\Mmod{2N}.
\end{gather*}
Observe that
\begin{equation*}\Fma\Fmn^{-1}=AN^{-1}\BmZ+\frac{B+\sqrt{D}}{2}\BmZ~;~
N\tau_\Fmn^\SmA=\tau^{\SmA[\Fmn]^{-1}} \Mmod{\Gamma}.
\end{equation*}

%\subsubsection*{Integral quadratic forms}
From either of the previous description we extract the following relations\footnote{there is a small unimportant misprint in \cite[p.~236]{GZ}.} between Heegner points of level $1$ and level $N$
\begin{equation}\label{HeegnerLevel}
\sigma\cdot \tau_\Fmn^{\SmA}=\tau_{\overline{\Fmn}}^{\SmA[\Fmn]^{-1}}~;~
\iota(\tau_\Fmn^\SmA)=(\tau^\SmA,\tau^{\SmA[\Fmn]^{-1}}).
\end{equation}

\subsection{Changing levels.}
The expression \eqref{startn} does not depend on the choice of
$\Fmn$. Hence we could without loss of generality enlarge the average to the $\Fmn$. We
avoid this for two reasons. This is useless and in general Waldspurger formula and its variants involve a toric orbit in the adelic set-up that corresponds in our classical set up to a \emph{single orbit} over ${\MCl}_K$.

By \eqref{HeegnerLevel}, we have:
\begin{equation}
E(s,\tau^{\SmA[\Fmn]})=E(s,N\cdot\tau_{\overline{\Fmn}}^{\SmA})
\end{equation}

Hence we introduce the function
\begin{equation*}
\psi(s,z):=|\zeta(2s)|^2 E(s,z)\overline{E(s,Nz)}=\zeta(s)\zeta(1-s)E(s,z)E(1-s,Nz)
\end{equation*}
in the variable $z\in Y_0(N)$ so that we need now to evaluate the ``period``:
\begin{equation}
S:=\frac{1}{h}\suma{\psi(s,\tau^\SmA_\Fmn)}.
\end{equation}

\subsection{Sketch of proof.}\label{subsec:sketch}
We introduce a Hecke operator which explains where the two terms $N^{-\Mdemi}\log N$ and $N^A|D|^{-\delta}$ are coming from. In this paragraph we shall ignore that the function $\psi$ is not $L^1$ and give a brief sketch. In the next subsection we produce a rigorous argument.

Duke's Theorem would yield (if $\psi$ were smooth):
\begin{equation*}
S = \frac{3}{\pi(N+1)}\int_{Y_0(N)}{\psi(z)\frac{dxdy}{y^2}}+O(N^A D^{-\delta}),
\end{equation*}
the $N^A$ coming from uniform bounds on the spectral coefficients of
$\psi$.

Observe that, at least formally, by definition of the Hecke correspondence, we have:
\begin{equation*}
\int_{Y_0(N)}{\psi(z)\frac{dxdy}{y^2}}=\int_{Y(1)}{E(s,z)\overline{\TmT_N\cdot E(s,z)}\frac{dxdy}{y^2}}.
\end{equation*}
The Eisenstein series $E(s,z)$ is an eigenvalue of the Hecke operator
$\TmT_N$ with eigenvalue $N^s+N^{1-s}$. Since $\Re s=\Mdemi$, we are
left with a factor $N^{-\Mdemi}$ times the integral
of $\left|E(s,z)\right|^2$ which we do as if it were finite.
We finally would obtain the bound
$S \ll N^{-\Mdemi}+ N^A D^{-\delta},$
we were looking for.

\subsection{Regularization s.}
We introduce the function $R(s_1,s_2,z)$ defined in the 2-dimensional strip $\frac{1}{4}<\Re s_1,\Re s_2<\frac{3}{4}$ by the following equality (observe that the factor $\pi^{-s_1-s_2}\Gamma(s_1)\Gamma(s_2)$ does not vanish in that strip):
\begin{equation}\begin{aligned}\label{regularized}
\pi^{-s_1-s_2}\Gamma(s_1)\Gamma(s_2)
&R(s_1,s_2,z):=\\
  &\xi(2s_1)\xi(2s_2)[N^{s_2}E_\infty(s_1+s_2,z)+N^{s_1}E_0(s_1+s_2,z)]+\\
  +&\xi(2-2s_1)\xi(2s_2)[N^{s_2}E_\infty(1+s_2-s_1,z)+N^{1-s_1}E_0(1+s_2-s_1,z)]+\\
  +&\xi(2s_1)\xi(2-2s_2)[N^{1-s_2}E_\infty(1+s_1-s_2,z)+N^{s_1}E_0(1+s_2-s_1,z)]+\\
  +&\xi(2-2s_1)\xi(2-2s_2)[N^{1-s_2}E_\infty(2-s_1-s_2,z)+N^{1-s_1}E_0(2-s_1-s_2,z)].
\end{aligned}\end{equation}
It is possible to check by various means that it is actually holomorphic in the $s_1$ and $s_2$ variables. Here is a possibility. The potential singularities are located on the union of the four complex lines
\begin{equation*}
  \Bigl(s_1=\Mdemi \Bigr) \text{ or } \Bigl(s_1+s_2=1 \Bigr) \text{ or } \Bigl(s_2=\Mdemi \Bigr) \text{ or } \Bigl(s_1=s_2 \Bigr)
\end{equation*}
and they are at most simple. One checks by computing the residues or better thanks to the symmetries ($s_1\leftrightarrow 1-s_1$) and ($s_2\leftrightarrow 1-s_2$) that $R$ has no singularity on those lines, except perhaps at $(s_1,s_2)=(\Mdemi,\Mdemi)$. But this potential singularity isolated at $(\Mdemi,\Mdemi)$ is a removable singularity by Hartog's Lemma.

Then we regularize the function $\psi$ by forming the difference:
\begin{equation}\label{diff2}
  C(s,z):=|\zeta(2s)|^2 E(s,z)E(1-s,Nz)-R(s,1-s,z).
\end{equation}
This process cuts the sum $S$ into two parts. We first claim that:
\begin{equation}\label{temp:R}
  \frac{1}{h}\suma R(s,1-s,\tau_{\Fmn}^{\SmA}) \ll_\epsilon N^{-\Mdemi}(\log N)^3 |Ds|^\epsilon.
\end{equation}
Indeed, it is easy to evaluate $\suma R(s_1,s_2,\tau_\Fmn^\SmA)$ explicitly when $s_1$ and $s_2$ are generic. From the Hecke formula \eqref{Hecke} together with \eqref{eisensteinniveaux} we obtain for example that
\begin{align*}
  \frac{1}{h}\suma \zeta(2s)E_\infty(s,\tau_\Fmn^\SmA)
  &= N^{-s}(1+N^{-s})^{-1} 
  \frac{1}{2h}(\frac{\sqrt{|D|}}{2})^s\zeta(s)L(s,\chi_D) \\
  &\ll_\epsilon N^{-\Mdemi}\log N \cdot |D|^\epsilon \qquad \text{for $\Re s=1$ and $s\not=1$.}
\end{align*}
When we instead specialize the average of \eqref{regularized} at $s_1=s$ and $s_2=1-s$ it is not difficult to obtain \eqref{temp:R} (the exponent $3$ in the logarithms appear only at the point $s=1$).

We also claim that for some absolute constants $A,B,\delta>0$ we have
\begin{equation}
  \frac{1}{h}\suma C(s,\tau_{\Fmn}^{\SmA}) \ll |D|^{-\delta} N^A |s|^B
\end{equation}
thus finishing the bound for $S$ and concluding the proof of Theorem~\ref{th:dfiB}.

This claim is a consequence of Proposition~\ref{maintool} whose assumptions are fulfilled by the next Proposition. We leave the details of the proof to the reader since it is very similar to Proposition~\ref{prop:B}.

\begin{proposition}
The function $z\mapsto C(s,z)$ satisfies the following convex bounds on $N,\Re s=\Mdemi$ and $z\in Y_0(N)$. Here $y$ stands for the 'height' of $z$, that is $\max(\Mim z,\Mim \sigma z)$ in the standard fundamental domain.
\begin{align}
\tag{C1}
\frac{\partial^{i+j}}{\partial x^i\partial
  y^j}C(s,z) &\ll_{i,j} y^{A_{ij}} N^{B_{ij}} |s|^{C_{ij}}\quad\text{for some $A_{ij},B_{ij}, C_{ij}>0$};\\
\tag{C2}
C(s,z),\ C(s,\sigma\cdot z) &\ll N^A|s|^B \log y\quad\text{for some $A,B>0$};\\
\label{integrale}\tag{C3}
\int_{Y_0(N)}{C(s,z)\mu(dz)} &=0.
\end{align}
\end{proposition}

%--------------------------
\section{Applications.}\label{discussion}
%---------------------------
Obtaining polynomial bounds for ``periods`` in the spectral parameter is of central significance for applications. This is typically the case in the theory of shifted convolution sums (see \cite{Sarn94}). Here our applications are much more modest of course.

Consider for instance the restriction to the diagonal of general automorphic kernels:
\begin{equation}\label{automorphe}
K(z,z):=\sum\limits_{\gamma\in PSL_2(\BmZ)}{k(d(z,\gamma z))},
\end{equation}
where $k\in \CmC^\infty_0(\BmR^\times_+)$ is a smooth function of compact support and $d()$ is the hyperbolic distance. Spectral expansions of automorphic kernels appear in many different contexts. The continuous part reads:
\begin{equation*}
 \int_{\MRe s=\Mdemi} \widehat{k}(s) \abs{E(s,z)}^2 \frac{ds}{4i\pi}.
\end{equation*}
 The spectral coefficient satisfies: $\widehat{k}(s) \ll_A |s|^{-A}$. 
 
When evaluating the average over Heegner points $\tau^\SmA$, it is natural to plug in the asymptotic given in Theorem~\ref{th:dfiA}. To stay on the safe side, the polynomial growth of the error term in Theorem~\ref{th:dfiA} is clearly essential.

The following are two asymptotics taken from Chapter~10 of the author's PhD thesis \cite{thes:Temp}:
\begin{align}
\frac{1}{h(D)}\sum\limits_{\SmA\in\MCl_K}{\log |\Mim \tau^\SmA\cdot j'(\tau^\SmA)|}
&=6 \CmL_D + a_1 +O(|D|^{-\delta})\\
\frac{1}{h(D)}\sum\limits_{\SmA\in \MCl_K}{K(\tau^\SmA,\tau^\SmA)}
&= a_2 \CmL_D + a_3+O_k(|D|^{-\delta})
\end{align}
These results are to be compared with the exact average~\eqref{Main1} which is equal to $\CmL_D+\log 2-\gamma$. In all three cases the averaged function is not $L^1$ but barely and the asymptotics have similar shapes.

We made explicit the constants $a_1$, $a_2$ and $a_3$ in \cite{thes:Temp}. We have $a_2
=2\int_0^\infty k(u)u^{-\Mdemi}du\\
=\int_{(\Mdemi)}\widehat{k}(v)\frac{dv}{2i\pi}$, but the exact computation of the constants $a_1$ and $a_3$ is long and tedious. These asymptotics emerge in the analytic study of the Gross-Zagier formula \cite{GZ1,GZ}. Indeed the quantity $\log |\Mim z\cdot j'(z)|$ is closely related to a {\em regularized} height on $X_0(1)_\BmQ \simeq {\bf P^1_\BmQ}$. We shall return to this question in future papers.

%\nocite{*}
%\addcontentsline{toc}{section}{Bibliographie}
\bibliographystyle{amsplain}
%\bibliography{bibli.bib}
\providecommand{\bysame}{\leavevmode\hbox to3em{\hrulefill}\thinspace}
\providecommand{\MR}{\relax\ifhmode\unskip\space\fi MR }
\providecommand{\MRhref}[2]{%
  \href{http://www.ams.org/mathscinet-getitem?mr=#1}{#2}
}
\providecommand{\href}[2]{#2}

\end{document}